\documentclass[12pt]{article}
\usepackage{amsmath,amsfonts,latexsym}

\usepackage{amsfonts}

\newcommand{\C}{{\mathbb C}}

\newcommand{\Z}{{\mathbb Z}}

\newcommand{\G}{{\cal G}}

\begin{document}
\begin{center}
{\LARGE{\bf Partial classification of modules for Lie-algebra of 
diffeomorphisms of $d$-dimensional torus}}  \\ [5mm] 
{\bf S. Eswara Rao}\\
School of Mathematics \\
Tata Institute of Fundamental Research \\
Homi Bhabha Road \\
Mumbai 400 005 \\
India \\ [1.5cm]
email: senapati@math.tifr.res.in
\end{center}
\begin{abstract}

We consider the Lie-algebra of the group of diffeomorphisms of a $d$-
dimensional torus which is also known to be the algebra of derivations 
on a Laurent polynomial ring $A$ in $d$ commuting variables denoted by
${\rm Der}A$.  The universal central extension of Der $A$  for $d=1$ is
the so called Virasoro algebra. The connection between Virasoro algebra
and physics is well known.  See for example the book on Conformal Field 
Theory by Di Francesco, Mathieu and Senechal.

In this paper we classify ($A$, Der $A)$ modules which are irreducible and
has finite dimensional weight spaces.  Earlier Larsson constructed a large
class of modules the so called tensor fields based on $g \ell_d$ modules
which are also $A$ modules.  We prove that they exhaust all $(A$, Der $A)$
irreducible modules.

\end{abstract}

\vskip 1cm
\section*{Introduction}

\quad

It is well known that the group of  diffeomorphisms on a manifold 
is very important and shows up directly in many branches of physics
(see for example Ref.22).  We are particularly interested in
$d$-dimensional torus.  The case $d=1$ is well studied by both
mathematicians and physicists.  The one dimensional central extension of
the
Lie-algebra of diffeomorphisms of the circle is well known object called
Virasoro algebra.  The representation theory of Virasoro algebra is
studied in great detail.  See Ref.13.

The Virasoro algebra acts on any (except when the level is negative of
dual Coxeter number) highest weight module of the affine Lie algebra
through the use of the famous Sugawara operators.  It is well known that
affine Lie algebras admit representation on the Fock space (see Ref.12) and
hence admits a representation of the Virasoro algebra.  This classical
theory is what we originally want to generalize to $d$-dimensional torus.

The relation to physics is well established in the book on Conformal Field
Theory by Di Francesco, Mathieu and Senechal. Ref.10.  Specially the
chapters 13 to 18 explains the connection between
physics and the representation theory of Virasoro and affine Kac-Moody
Lie-algebras.  Several important papers on these aspects have been put in one
volume by Goddard and Olive in Ref.8.  The most fundamental paper in this
direction is due to Belavin, Polyakov and Zamolodchikov in Ref.1.

The generalization of affine Lie algebra is the  so called toroidal
Lie-algebra.  For the first time a large class of representations are 
constructed in Refs. 6 and 20  through the use of vertex operators
generalizing the Fock space construction to the toroidal Lie-algebras.
One significant difference for the toroidal Lie-algebra is that the
universal center is infinite dimensional unlike in the affine case where
it is one dimensional.

So the next natural question is to generalise Virasoro algebra and see
whether the algebra acts on the Fock space.  For that we first denote
the Lie algebra of diffeomorphisms of $d$-dimensional torus by DerA (it
is known that DerA is isomorphic to the derivations of Laurent polynomial
ring A in d-variables).  Here one should mention that several 
attempts  has been made by physicists to give a Fock space representation
to DerA or to its extension (see Ref. 9).  They all failed to produce
any interesting results due to lack of proper definition  of ``normal
ordering'' among other things.  At this juncture an interesting result has
come out from Ref.21 which says that DerA has no non-trivial central
extension for $ d \geq 2$.

Let us go back to the vertex construction of toroidal Lie-algebra of Ref.
6. Here operators are constructed for DerA generalizing the Sugawara
Construction. But  the corresponding extension for DerA is very wild
(certainly non-central) and not tractable (see Ref.5). In the process an
interesting abelian extension for DerA has been created in Ref. 6 and the
abelian part is exactly the center of the toroidal Lie algebra.  So the
semi-direct sum of the toroidal Lie algebra and DerA with common
extension has emerged as an interesting object which we will now define.

We will first define toroidal Lie algebra.  Let $\G$ be simple finite
dimensional Lie algebra and let $<, >$ be a non-degenerate symmetric
bilinear form on $\G$.  Fix a positive integer d and let $A = \C
[t_1^{\pm 1}, \cdots, t_d^{\pm 1}]$ be a Laurent polynomial ring in
d-commuting variables.

Let $\Omega_A$ be the module of differentials which can be defined as
 vector space spanned by $t^r K_i, \ i=1, \cdots, d$ and $r \in \Z^n$.  Let
$d_A$ be the subspace spanned by $\Sigma r_i t^r K_i$ and consider the 
toroidal Lie algebra $\G \otimes A \oplus \Omega_A/d_A$ with Lie-bracket.
$$
 [X \otimes t^r, Y \otimes t^s]=[X,Y] \otimes t^{r+s}+ <X,Y>  \Sigma r_i
t^{r+s} K_i$$
$ \Omega_A/d_A$ is central.

Let DerA be the derivation on A. For $ u \in \C^d$ and $r \in \Z^d$ 
let
$ D(u,r) = \Sigma u_i t^r t_i \frac{d}{dti}$ where
$u=(u_1,  \cdots, u_d)$ and $r =(r_1, \cdots, r_d), t^r=t^{r_1}_1 t_2^{r_2}
\cdots  t_d^{r_d} \in A$. Let $K(u,r) = \Sigma u_i t^r K_i$.  Consider the
Lie-algebra
$$
\tau = \G \otimes A \oplus \Omega_A/d_A \oplus DerA
$$ 
$$
[D(u,r), D(v,s)] = D(w, r+s)-(u,s)(v,r) K (r, r+s)$$
where
$w=(u,s)v -(v,r)u$ and ( , ) is the standard inner product in $\C^d $
$$[D(u,r), K(v, s)]= (u,s)K(v,r + s)+(u,v)K(r,r+s)
$$
$$
 [K(u,r),K(v,s)]=0
$$
The first question is that can we construct a representation for $\tau$ from 
known methods. Several attempts have been made (Ref. 2, 3 and 6).  
Eventually in a remarkable paper Yuly Billig (Ref.24) has succeeded in 
constructing a class of modules for $\tau$ through the use of Vertex 
operator algebras (VOA).  In the process Yuly Billig has used the DerA 
modules constructed in Ref.16 and studied in Ref.4. One natural question
is, 
does there exists more modules for DerA so that we get a much larger 
class of models of $\tau$.

In an interesting paper by Jiang and Meng (Ref.11) it is proved that 
classification of irreducible integrable modules of $\tau$ can be reduced 
to the classification of irreducible (A, DerA) modules which the current 
paper settles.  See also ref.7 for more precise results.

Let me explain the results of this paper in more detail. In Ref.16, 
Larsson constructed a large class of DerA modules and some of them with 
finite dimensional weight spaces.  In fact he constructed a functor from 
$gl_d$-modules to DerA-modules. In Ref.4 the author proves that the 
image of an irreducible finite dimensional module is most often irreducible.

Further Larsson's DerA modules are A modules too and they are always
irreducible as (A, DerA)-modules.  Thus the purpose of the paper is to
prove the converse.   So we prove in Theorem 1.9 that any (A, DerA)
module which is irreducible and has finite  dimensional weight  spaces has
to come from Larsson's construction.

It will be certainly interesting to classify  all irreducible DerA
modules with finite dimensional weight spaces.   Now some kind of highest
weight modules are constructed in  Ref.2 and 3  (They are not A modules).
We will also note that $GL(d, \Z)$ acts as  automorphisms on DerA and so
we can twist a DerA  module  by $GL(d,\Z)$. (Larsson's modules are closed
under $GL(d,\Z)$ action).  So  it will be interesting to prove that any
irreducible DerA module with  finite dimensional weight spaces is either
an highest weight module or a Larsson's module upto a twist of  
$GL(d,\Z)$ action.  The problem is completely solved for d=1 by Mathieu
in Ref.19

Our  results in Ref.6 have been reinterpreted in the language of physics
by  Larsson in Refs. 17 and 18. Larsson is also first to talk about 
non-central extensions is Ref. 15. DerA has also been studied in Ref.23.

\section*{1. Section}

\paragraph*{(1.1)}  Throughout this paper we fix a positive integer $d
\geq 2$ and a Laurent polynomial ring $A=\C[t_1^{\pm 1}, \cdots t_d^{\pm
1}]$ in $d$ commuting variables.  Let $\C^d$ be $d$ copies of complex
field $\C$.  Let $e_1, \cdots e_d$ be the standard basis of $\C^d$ and let
$(,)$ be the standard form on $\C^d$ such that $(e_i, e_j)=\delta_{ij}$.

\paragraph*{(1.2)}  Let $\Gamma =\Z e_1 \oplus \cdots \oplus \Z e_d$.  
Throughout this paper we use $m,n,r$ and $s$ to denote elements of 
$\Gamma$.  For $r=\sum r_i e_i \in \Gamma$ let $t^r =t_1^{r_1} t_2^{r_2}
\cdots t_d^{r_d} \in A$ and  let $D^i (r) = t^r t_i  \frac{d}{dt_i}$ be
 a derivation on $A$.  Let ${\rm Der}A$ be the Lie-algebra of derivations
of $A$.  It is easy to verify that $D^i (r), 1 \leq i \leq d, r \in \Gamma$
is a basis of ${\rm Der}A$.  For $u= \sum u_i e_i \in \C^d$
let $D(u,r) = \sum u_i D^i (r)$.  Then ${\rm Der}A$ has the following Lie
structure:

\paragraph*{(1.3)}
$$[D (u,r), D(v,s)] =D(w, r+s)$$
where $w= (u,s) v - (v,r) u, \ r,s \in \Gamma$ and $u,v \in \C^d$.  Let
$h$ be the subspace spanned by $D^i (0), 1 \leq i \leq d$, which is a 
maximal abelian subalgebra of ${\rm Der}A$.

\paragraph*{(1.4)}  Note that $D(u,r) t^s = (u,s) t^{r+s}$.
Thus $A \oplus {\rm Der}A$ is a Lie-algebra by extending the Lie structure
in the following way
$$[t^r, t^s]=0$$
$$[D(u,r),t^m] =(u,m) t^{r+m}.$$
Let $\tilde{h} =\C \oplus h$ which is an abelian subalgebra of $A \oplus {\rm
Der}A$.

The purpose of this paper is to study $A \oplus {\rm Der}A$ modules which
are weight modules for $\tilde{h}$ with finite dimensional weight spaces
and to classify such modules with some natural conditions.

We first recall ${\rm Der}A$ modules which  are constructed and studied in 
References 16 and 4.

\paragraph*{(1.5)}  Let $g \ell_d$ be the Lie-algebra of $d \times d$ matrices
with entries in $\C$.  Let $E_{ij}$ be the elementary matrix with $(i,j)$
th entry 1 and zero elsewhere.  Then it is well known that $g \ell_d$ is
spanned by $E_{ij}, 1 \leq i, j \leq d$ with the following Lie-bracket
$$[E_{ij}, E_{k \ell}]= \delta_{jk} e_{i \ell} - \delta_{i \ell} E_{kj}.$$
Let $g \ell_d = s \ell_d \oplus \C I$ where $s  
\ell_d$ is a Lie-subalgebra of trace zero matrices and $I$ is the
identity matrix.  Let $V (\psi)$ be the irreducible finite dimensional
module for $s \ell_d$
where $\psi$ is a dominant integral weight.  Let $I$ act by  scalar
$b$ on $V(\psi)$ and denote the resultant $g \ell_d$ module by $V(\psi,b)$.  
Let  $\alpha \in \C^d$ and we will make $F^{\alpha} (\psi, b):= V (\psi, b)
\otimes A$ a ${\rm Der}A$ module.  First denote $v \otimes t^m$ by $v (m)$
for $v$ in $V(\psi, b)$ and $m$ in $\Gamma$.

\paragraph*{(1.6) Definition  } (Ref. 16)
$$D(u,r) \cdot v (m) = (u, m+\alpha) v (m+r)+ (\displaystyle{\sum_{i,j}}
u_i r_j E_{ji} v) (m+r)$$
where $m,r \in \Gamma, u \in \C^d, v \in V(\psi,b)$.  We will now recall
the following

\paragraph* {(1.7) Theorem}  (Theorem (1.9) and Proposition  (5.1) of
Ref.4). \\
(1) $F^{\alpha} (\psi, b)$ is irreducible as ${\rm Der}A$ module if
$(\psi, b) \neq (\delta_k, k), (0,b), 1 \leq k \leq d-1$ where $\delta_k$
is the $k$th fundamental weight of $s \ell_d$. \\
(2) $F^{\alpha} (0,b)$ is irreducible as ${\rm Der}A$ module unless $\alpha
\in \Gamma$ and $b \in \{0,d\}$.  

In all other cases $F^{\alpha} (\psi,b)$ is reducible and the submodule
structure has been worked out in Proposition (5.1) and Theorem (5.5) of 
Ref.4.

Recall that $A$ is associative algebra with unit and $F^{\alpha} (\psi,b)$ 
is a $A$ module by defining
$$t^m \cdot v (r) =v (m+r)$$
for $m,r \in \Gamma$ and $v \in V(\psi,b)$.  Further it is easy to see that
 $F^{\alpha} (\psi, b)$ is a $A \oplus {\rm Der}A$ module.

\paragraph*{(1.8) Proposition}  $F^{\alpha} (\psi, b)$ is irreducible as
$A \oplus {\rm Der}A$ module.

\paragraph*{Proof} First note that $F^{\alpha} (\psi, b)$ is a 
weight module with respective to $\tilde{h}$ and the weight spaces are 
$V(\psi) \otimes t^m$.  Suppose $W$ is a non-zero $A \oplus {\rm Der}A$
submodule of $F^{\alpha}(\psi, b)$.  As submodule of a weight module is a  
weight module, $W$ is a  weight module.  From the action of $A$ it is 
clear that $v (m) \in W$ implies $v (s) \in W$ for all $s \in \Gamma$.  
Thus $W=W_1 \otimes A$ for some $W_1 \subseteq V (\psi,b)$.  Now choose
$u=e_i, r=e_j $ and 
consider
$$D (u,r) v (m) =(u,m+\alpha) v (m+r) + (E_{ji} v) (m+r).$$
It now follows from the above remarks that $W_1$ is $g \ell_d$-invariant.
Since $V(\psi)$ is irreducible and $W_1$ is non-zero, it follows that
$W_1=V (\psi)$ and  hence $W=F^{\alpha} (\psi,b)$.

The purpose of this paper is to prove converse of the above proposition.
In other words we classify $A \oplus {\rm Der}A$ modules with certain 
natural properties.

\paragraph*{(1.9) Theorem}  Let $V$ be irreducible module for $A \oplus 
{\rm Der}A$ which is also a weight module for $\tilde{h}$ with finite
 dimensional weight spaces.  We further assume the following:

\begin{enumerate}
\item[(1)] $V$ is a $A$-module as associative algebra and the Lie-module
structure of $A$ comes from associative algebra.

\item[(2)] $1. v =v, \forall v$ in $V$. \\
  
Then $V \cong F^{\alpha} (\psi,
b)$ (for some $\alpha, \psi, b)$ as $A \oplus {\rm Der}A$-module.
\end{enumerate}

We need to develop several lemmas to prove the theorem which will be done
in Section 2.  The final proof will be given in Section 3.

\section*{2. Section}
\quad \quad
First we need to change some notation.  We treat $A$ as group algebra over 
$\Gamma$.  For that let $k (r)$ be a symbol for $r \in \Gamma $.  Let
$A$ be the linear span of $k (r), r \in \Gamma$ with multiplication 
defined as $k (r) \cdot k (s) = k (r+s)$.

Let $U$ be the universal enveloping algebra of $A \oplus {\rm Der}A$. 
Let $L$ be the two sided ideal of $U$ generated by $k(r) k (s)-k(r+s)$ 
and $k (0) -1$. 

Throughout this section the module $V$ is as in Theorem (1.9).  Since
$V$ is a $A$ module, $L$ acts trivially and hence $V$ is a $U/L$-module.
Let $V=\displaystyle{\oplus_{r \in \Gamma}} V_r$ be the weight space
decomposition and $V_r =\{v \in V \mid D(u, 0)v=(u,r+\alpha)v, \forall u$
in $\C^d \}$.  Such a uniform $\alpha$ in $\C^d$ exists as $V$ is
irreducible.
In fact take any weight space where $h$ acts as linear function which can
be taken as $u \mapsto (u,\alpha)$ for some $\alpha \in \C^d$. Because of
irreducibility the action of $h$ on the rest of the spaces is easily computed.
  Further each $V_r$ is a $\tilde{h}$-module as $1$ in $A$ acts as one on the
entire module.

\paragraph*{(2.1)}  Let $U_1 = U/L$ and let $T(u,r)=k(-r) D(u,r)-D(u,0)$
as an element of $U_1 $ for $u \in \C^d$ and $r \in \Gamma$.  Let $T$ be
the subspace spanned by $T(u,r)$ for all $u$ and $r$.

\paragraph*{(2.2) Proposition}
\begin{enumerate}
\item[(1)]  $[T (v,s), T (u,r)]=(u,s) T(v,s)-(v,r)T(u,r)+T(w,r+s)$ \\
where $w=(v,r)u- (u,s)v$ and hence $T$ is a Lie-subalgebra.
\item[(2)] $[D (v,0), T(u,r)]=0$
\item[(3)] $V = \oplus V_r$ be weight space decomposition.  Then each
$V_r$ is $T$-invariant
\item[(4)] Each $V_r$ is $T$-irreducible
\item[(5)] $V_r \cong V_s$ as $T$-module.
\end{enumerate}

\paragraph*{Proof}
\paragraph*{(2)} $[D (v,0),k (-r) D(u,r)-D(u,0)]$\\
$$
\begin{array}{lll}
&=&[D (v,0), k(-r)D (u,r)] \\
&=&[D (v,0), k(-r)]D (u,r) \\
&+&k (-r)[D (v,0),D(u,r)]\\
&=&-(v,r) k (-r)D (u,r) \\
&+& (v,r) k (-r) D (u,r)\\
&=& 0.
\end{array}
$$
(1) From (2) it follows that
$$[T (v,s), T(u,r)]=[k (-s) D(v,s), k (-r)D(u,r)]$$
$$=[k(-s), k(-r) D(u,r)]D (v,s)+k(-s) [D(v,s), k(-r)D (u,r)]$$

$$
\begin{array}{lll}
&=& [k(-s), k(-r)] D(u,r) D(v,s) \\
&+& k (-r) [k (-s),D(u,r)]D(v,s) \\
&+& k (-s) [D (v,s), k(-r)]D(u,r)\\
&+& k(-s) k (-r) [D (v,s),D(u,r)] \\
&=& (u,s) k(-s) D(v,s)\\
&&-(v,r) k (-r)D(u,r) \\
&+& k (-s-r) D(w,r+s) \ {\rm where} \ w=(v,r) u-(u,s) v \\
&=& (u,s) T(v,s) \\
&&-(v,r) T (u,r) \\
&+& T(w, r+s) \\
{\rm where} \ w&=& (v,r)u- (u,s) v
\end{array}
$$
(3) ~~ From (2) it follows that $T$ commutes with $h$ and hence $V_r$ is
a $T$-module. \\
(4) ~~ Let $U= \displaystyle{\oplus_{r \in \Gamma}} U_r$ where
$U_r =\{v \in U \mid [D (u,0), v] =(u,r)v$ for $u \in \C^d \}$.
Since $V$ is $A \oplus {\rm Der}A$ irreducible for $v,w$ in $V_r$ there
exists $X$ in $U_0$ such that $ X v=w$.  This is due to weight reasons.  Now
$X=\sum a_i X_i$ where each $X_i$ is of the form $k(-r) D(u_0,r_1) \cdots
D(u_k, r_k)$ where $\sum r_i=r$.  We are using the fact that $L$ acts
trivially on $V$.  Now using the fact that $k(-s) D(u,r)= D(u,r) k(s)-
(u,s) D(u,r)$ and the fact that $k (r) k (s) =k(r+s)$ we see that each $X_i$ is 
linear combination of elements of the form
$$k(-r_1) D(u_1,r_1) k (-r_2)D(u_2, r_2) \cdots k(-r_k)D(u_k,r_k).$$
This proves $X \in U(T)$, the universal enveloping algebra of $T$.  Hence
$V_r$ is  $T$ irreducible. \\
(5)~~ First note that $k (s-r) V_r \subseteq V_s$.  Repeating the same we
see that
$$V_r = k (r-s) k (s-r) V_r \subseteq k (r-s) V_s \subseteq V_r.$$
Thus $V_r = k (r-s) V_s$.  Define $f:V_r \to V_s$ by  $f(v) =k (s-r)v$ 
which is clearly injective and surjective.  Now
$$
\begin{array}{lll}
f (T (u,k)v) &=& k (s-r) T(u,k)v \\
&=& T (u,k) k (s-r)v \\
&=& T(u,k) f(v)
\end{array}
$$
Thus $f$ is a $T$-homomorphism.  This proves (5).

\paragraph*{(2.3) Notation}  For any integer $k>0, r,m_1, \cdots m_k \in 
\Gamma$ define 
$$
T_k (u,r, m_1, \cdots m_k)= T(u,r) - \displaystyle{\sum_{i}} T(u,r+m_i)
+ \displaystyle{\sum_{i <j}} T(u,r+m_i+m_j) \cdots 
$$
$$(-1)^j \displaystyle{\sum_{1 \leq i_1 < i_2 < \cdots < i_j \leq k}} T
(u,r+m_{i_1}+ \cdots + m_{i_j}) \cdots (-1)^k T(u,r+m_1+m_2+ \cdots
+m_k).$$
Let $I_k$ be the linear span of $T_k (u,r, m_1 \cdots m_k), u \in \C^d$ for
all $r, m_1, m_2 \cdots m_k \in \Gamma$.

\paragraph*{(2.4) Lemma}  

$(1) ~~ T_k (u,r,m_1, \cdots m_k) =T_k (u,r,m_{\sigma (1)}, \cdots 
m_{\sigma (k)})$
for any permutation $\sigma $ on $k$-letters. \\
(2) $T_k (u,r,m_1, \cdots m_k) = T_{k-1} (u,r, m_1, \cdots m_{k-1})
-T_{k-1} (u,r+m_k, m_1, \cdots m_{k-1})$ \\
(3) $I_k$ is an ideal of $T$. \\
(4) $I_k \subseteq I_{k-1} $ for $k \geq 2$. \\
(5) $[I_k, I_{\ell}] \subseteq I_{k+\ell-1} $ for $k, \ell \geq 1$.

\paragraph*{Proof} (1) Follows from definition \\ 
(2) Collect all terms where $m_k$ does not occur in the sum of $T_k$ and
that can be  seen to be equal to $T_{k-1} (u,r,m_1, \cdots m_{k-1})$.  Sum 
of the rest of the terms can be seen to equal to $-T_{k-1} (u,r+m_k,m_1, \cdots
m_{k-1})$.  This is because every term contains $m_k$ \\
(3) $[T(v,s),T_k (u,r,m_1, \cdots m_k)]=$ \\ 
$$(u,s) \displaystyle{\sum_{\ell=0}^{k}} 
\begin{pmatrix} k \\ \ell \end{pmatrix} (-1)^{\ell} 
k (-s) D(v,s)$$
$$ - (v,r) T_k (u,r,  m_1,  \cdots m_k)$$
$$+ \sum (v,m_i) T_{k-1} (u, r+m_i, m_1, \cdots \widehat{m_i}, \cdots m_k).$$
$$+(v,r) T_k (u, r+s,m_1, \cdots m_k)$$
$$-\sum (v,m_i)T_{k-1} (u,r+s+m_i,
m_1, \cdots \widehat{m_i}, \cdots m_k)$$
$$-(u,s)T_k (v,r+s,m_1, \cdots m_k).$$
By applying Proposition 2.2 (1) write $[T (v,s), T_k (u, r, m_1, \cdots
m_k)] =
A_1+A_2+A_3$.  It is easy to see that $A_1$ is the first term of the above 
formula.  Now in $A_2$ look for the terms where $(v,m_i)$ occurs and that 
can be seen as a component of the third term of the above formula.  Now
in $A_2$ the terms where no $m_i$ occurs is equal to the second term of
the above formula.  The rest of the formula can be seen in a similar way.  
This proves the claim.  Now note that the first term in the claim is zero. 
Clearly 2nd, 4th and 6th terms are in $I_k$.  Now 3rd and 5th term is
equal to
$$\sum (v,m_i) T_{k-1}(u, r+m_i, m_1, \cdots \widehat{m_i}, \cdots m_k)$$
$$-\sum (v,m_i) T_{k-1} (u, r+m_i+s, m_1 \cdots \widehat{m_i}, \cdots 
m_k)$$
$$=\sum (v,m_i) T_k (u,r+m_i, m_1, \cdots \widehat{m_i} \cdots m_k,s)$$
(by Lemma 2.4 (2)). \\
(4) Follows from (2). \\
(5)  
$$[T_\ell (v,s,n_1, \cdots n_{\ell}),T_k (u,r, m_1, \cdots m_k)]$$
$$\displaystyle{\sum_{t=0}^{\ell}} \ \  \displaystyle{\sum_{i_1 <i_2 <
\cdots <i_t}}(u,s+n_{i_1}+ \cdots + n_{i_t}).$$
$$\displaystyle{\sum_{b=0}^{k}} \begin{pmatrix}k \\b \end{pmatrix} (-1)^b
T(v,s+n_{i_1}
+ \cdots n_{i_t})$$
$$-\displaystyle{\sum_{t=0}^{k}} \ \displaystyle{\sum_{j_1 < \cdots < j_t}}
(v,r+m_{j_1}+ \cdots + m_{j_t}).$$
$$\displaystyle{\sum_{b=0}^{\ell}} \begin{pmatrix}\ell \\ b \end{pmatrix} 
(-1)^b T(u,r+m_{i_1}+ \cdots + m_{i_t})$$
$$+(v,r)T_{k+\ell} (u,r+s, m_1, \cdots m_k, n_1, \cdots n_{\ell})$$
$$-(u,s)T_{k+ \ell} (v,r+s,m_1, \cdots m_k,n_1, \cdots n_{\ell})$$
$$-\sum (v,m_i)T_{k+ \ell-1} (u,r+s+m_i, m_1, \cdots \widehat{m_i}, \cdots
m_k, n_1, \cdots n_{\ell})$$
$$+ \sum (u, n_j) T_{k+ \ell-1} (v,r+s+n_j, m_1, \cdots m_k, n_1,
\cdots \widehat{n_j}, \cdots n_{\ell})$$
The above formula can be deduced as in (3) from Proposition 2.2 (1).  Now
note that the first two terms are zero as
$$\displaystyle{\sum_{b=0}^{k}} \begin{pmatrix} k \\ b \end{pmatrix}
 (-1)^b =0 =
\displaystyle{\sum_{b=0}^{\ell}} \begin{pmatrix} \ell \\ b \end{pmatrix} 
(-1)^b.$$
The rest of the four terms are in $I_{k+\ell-1}$ and this proves (5).

\paragraph*{(2.5) Lemma}  For $u \in \C^d, 0 \neq m_i \in \Gamma, \ s \in 
\Gamma$. \\
(1)$T_k (u,s,m_1, \cdots m_k) \notin I_{k+1}$ for $k \geq 1$. \\
(2) $T_k (u,s,m_1, \cdots m_k)$
$$+ T_k (u,s,n, m_2, \cdots m_k)=T_k(u,s,m_1+ n, m_2, \cdots
m_k)+I_{k+1}$$
(3) $T_k (u,s, -m_1, m_2, \cdots m_k)$ \\
$$=-T_k (u, s-m_1, m_1, m_2, \cdots m_k)$$
\paragraph*{Proof}  To prove the Lemma, we first interpret $T_k$'s as
certain polynomials in $A=\C[t_1^{\pm 1} , \cdots t_d^{\pm 1}]$.  We fix a
non-zero $u$ in $\C^d$.  Let $k$ be a positive integer and let $m_1, m_2, 
\cdots m_k \in \Gamma$.  Let $P_k (m_1, \cdots m_k)=\displaystyle{\prod_{1 
\leq i \leq k}} (1-t^{m_i})$.  Recall $t^{m_i} = t_1^{(m_i)_1} \cdots
t_d^{(m_i)_d}$.
Let $J_k$ be the ideal in $A$ generated by $P_k (m_1, \cdots m_k)$ for
all non-zero $m_i$'s  $\in \Gamma$.  Then clearly $J_{k+1} \subseteq J_k$.
It is easy to see that  $T_k (u, r, m_1, \cdots m_k)$ can be identified with
polynomial $t^r P_k (m_1, \cdots m_k)$.  Recall that $u$ is fixed.

Thus it is sufficient to prove that,

\paragraph*{Claim 1} $P_k (m_1, \cdots m_k) \notin J_{k+1}$.  Suppose
$$P_k (m_1, \cdots m_k) = \sum f_{\ell} P_{k+1} (n_{\ell_1}, n_{\ell_2},
\cdots n_{\ell_{ k+1}}) \leqno{(*)}$$
where $f_{\ell} \in A$.
Let $Dt_i = t_i \frac{d}{dt_k}$.  Now consider $D_{t_{i_1}} \cdots
D_{t_{i_k}} P_{k+1} (n_1, \cdots n_{k+1})$ and evaluate at $(t_1, \cdots
t_d)= (1, \cdots 1)$.  This can be seen to be zero as after differentiating
 $P_{k+1}, k$ times, each component has at least one factor
$(t^{n_{i}}-1)$.
We will now prove that there exists $i_{1}, \cdots i_k$ such that

\paragraph*{Claim 2} $D_{t_{i_1}} \cdots D_{t_{i_k}} P_k (m_1,
\cdots m_k) \mid_{t=(1, \dots  1)}$ is non-zero.  Thus $*$ can not hold.  This
prove the claim 1.  Now choose $\ell, 1 \leq \ell \leq d$ such that
$S=\{i \mid (m_i)_{\ell} \neq 0\}$ is  non-empty.  Let $\# S=p$ and
let $i_1, \cdots i_p \in S$.  Consider
$$(D_{t_\ell})^p P_k (m_1, \cdots m_k)= \mu \displaystyle{\prod_{i \in S}}
(m_i)_{\ell} \displaystyle{\prod_{j \notin S}} (1-t^{m_j}) t^{m_{i_1}+ 
m_{i_2}+ \cdots + m_{i_p}}+J_{k-s+1}$$
which is not too difficult to see.  Where $\mu$ is a non-negative integer.
Repeating the process finitely many times (choosing different index 
$\ell^1 \neq \ell)$.  We see that there exists $i_1, \cdots, i_k 
$ such that $ D_{t_{i_1}}D_{t_{i_2}} \cdots D_{t_{i_k}}$.$P_k (m_1, \cdots, m_k)= \lambda t^{m_1+ \cdots
+m_k}+J_1$
where $\lambda$ is non-zero integer.  Now evaluating at $t=(1, \cdots 1)$
we see that claim 2 is true. 

To see (2) first note that
$$(1-t^m) (1-t^n)+ (1-t^{m+n})=(1-t^m)+ (1-t^n).$$
Thus
$$t^s \displaystyle{\prod_{i=2}^{k}} (1- t^{m_i}) (1-t^{m_1})+t^s
\displaystyle{\prod_{i=2}^{k}} (1-t^{m_i}) (1-t^n)$$
$$=t^s \displaystyle{\prod_{i=2}^{k}} (1-t^{m_i}) (1-t^n) (1-t^{m_1})$$
$$+ t^s \displaystyle{\prod_{i=2}^{k}} (1-t^{m_i}) (1- t^{m_1+n}).$$
This proves (2). \\
(3) is easy to check.
\paragraph*{(2.6) Lemma}
$$ {\rm Dim} (I_k / I_{k+1}) \leq d^{k+1}, k \geq 1 \leqno{(1)}$$
$$T =I_1 \leqno{(2)}.$$
In particular $I_k$ is a co-finite ideal in $T$.

\paragraph*{Proof}  First note that from Lemma 2.4(2) we have
$$T_k (u, r, m_1, \cdots m_k) = T_k (u,s,m_1 \cdots m_k) \ {\rm mod} \ 
I_{k+1}$$
 for all $r,s \in \Gamma$.  
Further
$$-T_k (u,0,m_1, \cdots m_k) =T_k (u,0,-m_1, m_2 \cdots m_k) \ {\rm mod} \
 I_{k+1}$$ 
which follows from above and Lemma 2.5 (3).  Now from additive
property of Lemma 2.5(2) it follows that $I_k/ I_{k+1}$ is spanned by
$T_k (u,0, e_{i_1}, \cdots e_{i_k})$ where $e_1, \cdots e_d$ is the
standard basis.  Thus (1) follows. (2) follows from definitions.  Now it
is easy to conclude that $I_k$ is a co-finite ideal for each $k$.

\section*{3. Section}

\quad \quad

We will explain the plan of the proof of Theorem (1.9).  First we will
prove that $T/I_2 \cong g \ell_d (\C)$.  Then we will prove that if
$I_k, k \geq 2$ is zero on a  finite dimensional irreducible module $V$ of
$T$ then $I_2$ is zero on $V$.  Thus $V$ is a module for $T/I_2
\cong g \ell_d (\C)$.  Further we prove that any co-finite ideal $J$ of
$T$
contains $I_k$ for large $k$.  Thus any irreducible finite dimensional module
$V$  of $T$ is actually a module for $T/I_2$.  From this it will be easy
to conclude Theorem 1.9 which will be explained at the end of the section.

 \paragraph*{(3.1) Proposition} $T/I_2 \cong g \ell_d (\C)$.  

\paragraph*{Proof} First recall
that $F^{\alpha} (\psi, b)$ is a $A \oplus {\rm Der}A$-module and each weight
space $V(\psi) \otimes t^m$ is a  $T$-module. It is easy to verify that $I_2$ 
acts trivially on $V(\psi) \otimes t^m$. Now note that $T(e_i, e_j) v(m)=
E_{ji} v (m) \neq 0$ for some $\psi$.  From this we conclude that
$T(e_i, e_j)$ is non-zero in $T/I_2$.  Now it is easy to see that
$T(u,s)+T(u,r) = T(u, r+s) \ {\rm mod} \  I_2$ and hence $T(e_i, e_j)$
spans $T/I_2$. Define $\pi: T/I_2 \to g \ell_d (\C)$.
$$\pi (T (e_i, e_j)) =E_{ji}$$
Consider
$$
\begin{array}{lll}
X&=& [T(e_i,e_j),T(e_k, e_{\ell})]\\
&=& [k (-e_j) D(e_i, e_j), k (-e_\ell)D (e_k, e_\ell)]\\
&=& -\delta_{i \ell} k (-e_{\ell}) D(e_k, e_\ell) \\
&&+ \delta_{kj} k (-e_j)D (e_i, e_j)\\
&-& \delta_{kj} k(-e_\ell -e_j) D(e_i, e_{\ell} + e_j)\\
&&+\delta_{i \ell} k (-e_\ell -e_j)D (e_k, e_\ell+e_j).
\end{array}
$$
Follows from Proposition (2.2).  Note that the following is true in $T/I_2$.
$$k(-e_\ell -e_j)D(e_s,e_\ell+e_j)=k(-e_{\ell}) D(e_s,e_{\ell})$$
$$+k(-e_j ) D(e_s, e_j) -D(e_s,0)$$
for $s=i,k$.  Thus $X= -\delta_{kj} (k (-e_\ell)D (e_i,e_\ell)-D(e_i,0))$ \\
$$+\delta_{i \ell} (k (-e_j)D (e_k, e_j)-D(e_k, 0))$$
$$= -\delta_{kj} T(e_i,e_\ell)+ \delta_{i \ell} T(e_k, e_j).$$
Thus $\pi$ defines a surjective homomorphism.   As $T(u,0)$ is zero it follows
that 
$T(e_i,e_j)$ span $T/I_2$ which proves dim $(T/I_2) \leq d^2$.  Thus $\pi$
defines an isomorphism.

\paragraph*{(3.2) ~~Lemma} (Yuly Billig)  Suppose ${\cal G}$ is a Lie-
algebra over  $\C$ and $J$ is an ideal with spanning set $J_{\alpha},
\alpha \in B$.  Suppose there exists an element $I$ in ${\cal G}$ such
that $[I, J_\alpha]=\lambda J_\alpha$, $\lambda \neq 0 $ for all $\alpha
\in
B$.  then $J$ acts trivially on any irreducible finite dimensional module
$V$
of ${\cal G}$.

\paragraph*{Proof} Since the base field is complex numbers and $V$ is
finite dimensional, $I$ has eigen vectors.  Let $\lambda_1, \cdots, \lambda_k$
be all the eigen values of $I$ on $V$.  Choose $\lambda_i$ such that
 $\lambda+\lambda_i$ is not an eigenvalue. Let $v$ be eigenvector with
eigenvalue $\lambda_i$ for $I$. Consider $I J_\alpha v=
J_{\alpha} Iv+[I, J_{\alpha}]v= (\lambda_i+\lambda) J_{\alpha} v$.  
This proves $J_{\alpha} v=0 \  \forall \alpha \in B$.  Let $W=\{w \in V 
\mid J_\alpha w=0 \ \forall \alpha \in B\}$.  Since $J$ is an ideal, it 
is easy to see that $W$ is a
${\cal G}$-module. But $W \neq 0$.  Since $V$ is irreducible $W=V$ which
proves that $J$ acts trivially on $V$.

 \paragraph*{(3.3) ~~ Proposition}  Suppose $V$ is irreducible finite
 dimensional module for $T$ such that $I_{k+1}$ acts trivially on $V$.
Then $I_2$ acts trivially on $V$.  

\paragraph*{Proof} From the proof of Lemma 2.4(3) we have
$$[T(v,s), T_k (u,r,m_1, \cdots m_k)=- (v,r)T_{k+1}(u,r, m_1, \cdots
m_k,s)$$
$$-(u,s)T_k (v, r+s, m_1, \cdots m_k) + \sum (v,m_i)T_k (u,r+m_i,m_1,
\cdots 
\widehat{m_i}, \cdots  m_k,s).$$
Let $I=\sum T(e_i,e_i)$ and note that $I$ is actually identity element
in $T/I_2  \cong g \ell_d (\C)$.  Thus $I$ is non-zero on $T/I_{k+1}$ for
$k \geq 1$.

\paragraph*{Claim} $$[I, T_k (u,r,m_1, \cdots m_k)=(k-1) T_k (u,r,m_1,
\cdots
m_k).$$
Consider 
$$[\displaystyle{\sum_{j}} T(e_j, e_j), T_k (u,r,m_1, \cdots m_k)
=-\sum u_j T_k (e_j, r+e_j, m_1, \cdots m_k)$$
$$+ \displaystyle{\sum_{i,j}} (m_i)_j T_k (u,r+m_i, m_1, \cdots \widehat{m_i}
, \cdots m_k, e_j).$$
Now we use Lemma 2.5(2) and the following facts.
\begin{enumerate}
\item[(1)] $T_k$ is linear in $u$
\item[(2)] $ T_k (u,r+m,m_1,  \cdots m_k) = T_k (u,r,m_1, \cdots m_k) \
{\rm mod} \ I_{k+1}$ (by lemma 2.4(2)).
\item[(3)] $I_{k+1}$ is zero on $V$.
\end{enumerate}
From that we conclude that
$$[I, T_k (u,r,m_1, \cdots m_k)]=-T_k (u,r,m_1, \cdots m_k)+kT_k(u,r,m_1,
\cdots m_k)$$
which proves the claim.  Now we can use Lemma (3.2) for the ideal 
$I_k$.  Thus $I_k$ is zero on $V$.   Repeating this argument we conclude
that $I_2$ acts trivially on $V$.  This argument breaks down for $k=1$ as we
cannot apply the Lemma 3.2.

\paragraph*{(3.4) ~~ Proposition}  Any co-finite ideal $J$ of $T$ contains 
$I_k$ for large $k$.

\paragraph*{Proof}  Claim $J \cap I_k$ is co-finite in $T$ for all $k$.  For 
that consider $\varphi:T \mapsto T/J \oplus T/I_k$
$$v \mapsto (v,v).$$
Clearly ker $\varphi = J \cap I_k $ and $T/{J \cap I_k}$  
is a subalgebra of finite dimensional Lie-algebra $T/J \oplus T/I_k$.  This
proves the claim.

Consider \\
$I_k / J \cap I_k \stackrel{\varphi_1}{\hookrightarrow} I_{k-1} / 
{J \cap I_k} \stackrel{\varphi_2}{\rightarrow} I_{k-1} / J \cap I_{k-1}$
where $\varphi_1$ is injective and $\varphi_2$ is surjective.  Let
$\rho_k =\varphi_2 \circ \varphi_1$.  Thus $\rho_k : I_k / {I_k \cap J}
\to I_{k-1}/ J \cap I_{k-1}$.  Clearly $\rho_k$ in injective.  Let
$t_n = {\rm dim} \ I_n/{I_n \cap J}$ and note that $t_{n+1} \leq t_n$.
Thus $\{t_n\}_{n \in \Z^+}$ is a decreasing sequence of non-negative
integers.  Therefore $t_k=s$ for some $s$ and for large $k>N$.

First we note the following two statements for a fixed $i$.

\begin{enumerate}
\item[(1)] For $\ell \neq i$ 
$$[T (e_{\ell}, - e_{\ell}), T_k (e_i, 0, e_{j_1}, \cdots
e_{j_k})]=-k_{\ell} T_k (e_i, 0, e_{j_1}, \cdots e_{j_k})$$
where $k_{\ell}$ is the number of $e_{\ell}$ that occur in $T_k (e_i, 0,
e_{j_1}, \cdots e_{j_k})$.

\item[(2)] Suppose the ideal $J$ contains $\sum a_{m,I}
T_m (e_i, 0, e_{j_1}, \cdots e_{j_m})$ where the number of $e_i$'s
that occur in $T_m (e_i, 0, e_{j_1}, \cdots e_{j_m})$ is same
for all $m$  where $I=\{j_1, \cdots j_k\}$.  Then $J$ contains $T_m (e_i,
0,e_{j_1}, \cdots e_{j_m})$ for $m \ni a_{m,I} \neq 0$.
\end{enumerate}
\paragraph*{(1)}  Follows from the proof of lemma 2.4(3).  (2) follows from 
(1).  We will prove the Proposition assuming $d \geq 3$ to avoid
some computations.  For a fixed $i$, consider the following set
$$S=\{T_k (e_i, 0, e_{j_1}, \cdots e_{j_k}) \mid  j_{\ell} \neq i \
{\rm for \ all} \ \ell \}.$$
Now choose $k \ni \# S>s$ and $k>N$.  Thus $S$ is linearly dependent mod
$I_k \cap J$.  Thus there exists non-zero scalars $a_I (I=\{j_1, \cdots
j_k\})$ such that
$$X= \sum a_I T_k (e_i, 0, e_{j_1}, \cdots e_{j_k}) \in J.$$
Now using (2) we conclude that
$$T_k (e_i, 0, e_{j_1}, \cdots e_{j_k}) \in J \ {\rm for \ some} \ I.$$
For $m,n \neq i$ consider
$$[T (e_m, e_n), T_k (e_i, 0, e_{j_1}, \cdots e_{j_k})] =
\ell \delta_{j_im} T_k (e_i, e_m, e_{j_1}, \cdots \widehat{e}_{j_i}, \cdots
e_{j_k}, e_n) \in J.$$
Now $T_k (e_i, e_m, e_{j_1}, \cdots \widehat{e}_{j_i}, \cdots e_{j_k}, e_n)=
T_k (e_i, 0, e_{j_1}, \cdots \widehat{e}_{j_i}, \cdots e_{j_k}, e_n)$
$$-T_{k+1} (e_i, 0, e_{j_1}, \cdots e_{j_i}, \cdots e_{j_k}, e_n).$$
Now by (2) it follows that
$$T_k (e_i, 0, e_{j_1}, \cdots \widehat{e}_{j_i}, \cdots e_{j_k}, e_n) \in J.$$
Now repeating this process we see that
$$T_k (e_i, 0, e_{j_1}, \cdots e_{j_k}) \in J \leqno{(*)}$$
for all possible indices $j_1, \cdots j_k$ which are all different from
$i$.  Applying $T (e_i, e_i)$ to the above vector to conclude
$$T_k (e_i, e_i, e_{j_1}, \cdots e_{j_k}) \in J$$
$$T_k (e_i, 0, e_{j_1}, \cdots e_{j_k}) - T_k ( e_i, e_i,e_{j_1}, \cdots
e_{j_k}) \leqno{(**)}$$
$$=-T_{k+1} (e_i, 0, e_{j_1}, \cdots e_{j_k}, e_i) \in J.$$
Fix $j \neq i$.  Replacing $k$ by $k+1$, consider the following vector which
is in $J$ by $(*)$.
$$[T (e_j, e_i), T_{k+1} (e_i, 0, e_{j_1}, \cdots, e_{j_{k+1}})] = P
\delta_{j j_{\ell}} T_{k+1} (e_i, e_j, e_{j_1}, \cdots \widehat{e}_{j_\ell},
\cdots e_{j_{k+1}}, e_i)$$
$$-T_{k+1} (e_j, e_i,e_{j_1}, \cdots e_{j_{k+1}}).$$
Now $T_{k+1} (e_i, e_j, e_{j_1}, \cdots \widehat{e}_{j_\ell}, \cdots
e_{j_{k+1}}, e_i)$ \\
$$= T_{k+1} (e_i, 0, e_{j_1}, \cdots \widehat{e}_{j_\ell}, \cdots
e_{j_{k+1}}, e_i) -T_{k+2} (e_i, 0, e_{j_1}, \cdots e_{j}, \cdots
e_{j_{k+1}},e_i).$$
Now by $(**)$ both vectors are in $J$.  Thus we conclude that 
$$T_{k+1} (e_j, e_i, e_{j_1}, \cdots e_{j_{k+1}}) \in J.$$
Now by (2) we see that
$$T_{k+1} (e_j, 0, e_{j_1}, \cdots e_{j_{k+1}}) \in J.$$
This is true for all possible indices $j_1, \cdots j_{k+1}$ which are
all different from $i$.

Now applying $T (e_{j_\ell}, e_i)$ for $j_\ell \neq j$ we see that
$$T_{k+1} (e_j, e_{j_\ell}, e_{j_1}, \cdots \widehat{e}_{j_\ell}, \cdots
e_{j_{k+1}}, e_i) \in J.$$

Now by  (2) we see that
$$T_{k+1} (e_j, 0, e_{j_1}, \cdots, \widehat{e}_{j_\ell}, \cdots e_{j_{k+1}},
e_i) \in J.$$
Repeating this process we see that $T_{k+1} (e_j, 0, e_{\ell_1}, \cdots
e_{\ell_{k+1}}) \in J$
for all possible $\ell_1, \cdots \ell_{k+1}$. 

Now using the technique  in the proof of Lemma 2.5  we see that
$$T_{k+1} (e_j, 0, m_1, \cdots m_{k+1}) \in J \ {\rm for \ all} \ m_i \in 
\Gamma $$
Now replasing $k+1$ by $k+2$ we see that
$$T_{k+2} (e_j, 0, m_1, \cdots m_{k+2})=T_{k+1} (e_j, 0, m_1, \cdots m_{k+1})
-T_{k+1} (e_j, m_{k+2}, m_1, \cdots, m_{k+1})\in J.$$
Then it follows that
$$T_{k+1} (e_j, m_{k+2}, m_1, \cdots m_{k+1}) \in J.$$
Strictly speaking we have it for non-negative coefficients. But the other
cases can be handled similarly. 
This proves $I_{k+1} \subseteq J$ and the Proposition.  Further
$s=0$.
\paragraph*{Proof of Theorem (1.9)}.   Let $V$ be a module as in Theorem.
Let $V = \displaystyle{\oplus_{r \in \Gamma}} V_r $ be the weight space
decomposition where
$$V_r =\{v \in V \mid D(u,0)v=(u, r+\alpha)v, \ \forall u \in \C^d\}.$$ 
We know that $V_r \cong V_s$ as $T$-modules from Proposition 2.2(5).
This  with Proposition (3.3), Proposition  (3.4) combined with the fact
that some co-finite ideal of $T$ acts trivially on $V_r$ tells us that
all $V_r'$s are isomorphic to some $V(\psi, b)$ as $g \ell_d$-modules.  
Note that the isomorphism between $V_r'$s  is given by $k(r)$ (from proof of 
proposition 2.2(5)).  Thus if we let $V_r = V(\psi, b) \otimes t^r$
we see that $k(r)v(s) = v (s+r)$ for $v $ in $V(\psi,b)$.

Now consider $T (u,r)$ in $T/ {I_2}$ and note that it is linear in both
variables.

Thus
$$
\begin{array}{lll}
T(u,r) v (s) &=& \displaystyle{\sum_{i,j}} u_i r_j T(e_i,e_j)v (s) \\
&=& \displaystyle{\sum_{i,j}} u_i r_j E_{ji} v (s)
\end{array}
$$
Therefore $k(-r) D(u,r) v(s) =D(u,0) v(s) +(\sum u_i v_j E_{ji}v) (s)$ \\
$$=(u, s+\alpha) v (s)+(\sum u_i r_j E_{ji} v) (s).$$
Multiply both sides by $k(r)$ we get
$$
\begin{array}{lll}
D(u,r)v (s) &=& (u, s+ \alpha) v (s+r) \\
&+& (\displaystyle{\sum_{ij}} u_i r_j E_{ji} v) (s+r).
\end{array}
$$
This completes the proof of the theorem.

\pagebreak

\begin{center}
{\bf REFERENCES}
\end{center}

\begin{enumerate}
\item A.A. Belavin, A.M. Polyakov and A.B. Zamolodchikov, Infinite
conformal symmetry in two dimensional quantum field theory, Nuclear
Physics B 241, 333-380 (1984).

 \item S. Berman, Y. Billig and J.
Szmigiclski, Vertex operator 
algebras and the representation theory of toroidal Lie-algebras, 
Contemporary Mathematics, Vol. 297 1-26 (2002).

\item S. Berman and Y. Billig, Irreducible  representations for toroidal
Lie-algebras, Journal of Algebra, 221, 188-231 (1999).

\item S. Eswara Rao, Irreducible representations of the Lie-algebra 
of the Diffeomorphisms of a $d$-dimensional torus, Journal of Algebra,
182, 401-421 (1996).

\item S. Eswara Rao, Generalized Virasoro operators (2004). To appear in 
Communications in Algebra.
\item S. Eswara Rao and R.V. Moody, Vertex  representations for n-toroidal
Lie algebras and a generalization of the  Virasoro algebra, Comm. Math.
Physics 159, 239-264 (1994).
\item S. Eswara Rao and C. Jiang, Classification of irreducible integrable
representations for the full toroidal Lie-algebras, preprint 2004.

\item P. Goddard and D. Olive, Kac-Moody and Virasoro algebras, World
Scientific, Singapore (1988).

\item F. Figueirido and E. Ramos, Fock space  representation of the algebra
of diffeomosphisms of the n-torus, Physics letters B, 238, 247-251 (1990).

\item P. Di Francesco, P. Mathieu and D. Senechal, Conformal Field Theory,
Springer Verlag, New York (1997).

\item C Jiang and D. Meng,  Integrable representations  for generalised
Virasoro-toroidal Lie-algebra, 270, 307-334 (2003).

\item V. Kac and I. Frenkel, Basic representation of  affine Lie algebras and
dual resonanace models, Invent. Math. 62, 23-66 (1980).

\item V. Kac and A.K. Raina, Bombay lectures on highest  weight
representations of infinite dimensional Lie-algebras,  World Scientific,
Singapore, (1987).

\item T.A. Larsson, Multi dimensional Virasoro algebra,  Physics Letters
B, 231, 94-96 (1989).

\item T.A. Larsson, Central and non-central extensions  of multi-graded
Lie algebras, Journal of Physics A, 25, 1177-1184 (1992).

\item T.A. Larsson, Conformal fields: A Class of  Representations of Vect
$(N)$, Internat. J. Modern Phys. A. Vol 7, No.26, 6493-6508 (1992).

\item T.A. Larsson, Lowest energy  representations  of non-centrally
extended diffeomorphism algebras, Comm. Math. Physics, 201, 461-470 (1999).

\item T.A. Larsson, Extended diffeomorphism algebras  and trajectories in
jet space, Comm. Math. Physics, 214, 469-491 (2000).

\item O. Mathieu, Classification of Harish-Chandra Modules over the
Virasoro Algebras, Invent. Math. 107, 225-234 (1992).

\item R.V. Moody, S. Eswara Rao and  T. Yokomuma, Toroidal  Lie algebras and
vertex representations, Geom. Ded. 35, 283-307 (1990).

\item  E. Ramos, C.H. Sah and R.E. Shrock, Algebra of diffeomorphisms
of the $N$-torus, J. Math. Phys.31, No.8., 1805-1816 (1990).

\item E. Ragoucy and P. Sorba, An attempt to relate area-preserving
diffeomorphisms to Kac-Moody Lie algebras, Letters in  Mathematical
Physics, 21, 329-342 (1991).

\item E. Ramos and R.E. Shrock, Infinite dimensional $\Z^n$-  indexed Lie
algebras and their super symmetric generalizations,  International Jour.
of Modern Phys. A, Vol.4, No. 16, 4295-4302 (1989).
\item Yuly Billig, Energy-momentum tensor for the  toroidal Lie algebras,
arXiv. Math. RT/0201313 (2002).
\end{enumerate}

\end{document}